%
%
%
\input amstex.tex
\documentstyle{amsppt}
\magnification=1200
\baselineskip=13pt
\NoBlackBoxes
\def\R{{\Bbb R}}

\def\C{{\Bbb C}}
\def\Z{{\Bbb Z}}
\def\Zp{{\Bbb Z}_+}
\def\a{\alpha}
\def\g{\gamma}
\def\as{\alpha^{\ast}}
\def\gs{\gamma^{\ast}}
\def\r{\rho_{\sigma,\infty}}
\def\rst{\rho_{\sigma,\tau}}
\def\hp{\hat p}
\def\hP{\hat P}
\def\P{P_l(\lambda;1,1,q^{2\sigma},1;q^2)}
\def\A{{\Cal A}_q}

\topmatter
\title Addition formula for big $q$-Legendre polynomials\\
from the quantum $SU(2)$ group\endtitle
\rightheadtext{Addition formula for big $q$-Legendre polynomials}
\author H.T. Koelink\endauthor
\affil  Katholieke Universiteit Leuven\endaffil
\address Department of Mathematics, Katholieke Universiteit Leuven,
Celestijnenlaan 200 B, B-3001 Leuven (Heverlee), Belgium\endaddress
\email erik\%twi\%wis\@cc3.KULeuven.ac.be\endemail
\thanks Supported by a NATO-Science Fellowship of the Netherlands
Organization for Scientific Research (NWO). \endthanks
\abstract From Koornwinder's interpretation of big $q$-Legendre
polynomials as spherical elements on the quantum $SU(2)$ group an
addition formula is derived for the big $q$-Legendre polynomial.
The formula involves Al-Salam--Carlitz polynomials,
little $q$-Jacobi polynomials and dual $q$-Krawtchouk polynomials.
For the little $q$-ultraspherical polynomials a product formula in terms
of a big $q$-Legendre polynomial follows by $q$-integration.
The addition and product formula for the Legendre polynomials
are obtained when $q$ tends to 1.
\endabstract
\endtopmatter
\document
\head 1. Introduction\endhead

Quantum groups provide a powerful approach to special functions of
basic hypergeometric type, cf. the survey papers by Koornwinder
\cite{12} and by Noumi \cite{16}, where the
reader will also find (more) references to the literature on
quantum groups and basic hypergeometric functions. In this paper
we show how the quantum group theoretic interpretation of
basic Jacobi polynomials leads to an addition formula
for the big $q$-Legendre polynomials
involving little $q$-Jacobi polynomials, dual $q$-Krawtchouk polynomials
and Al-Salam--Carlitz polynomials.

There are now several addition formulas available for basic analogues
of the Legendre polynomial. The addition formula for the
continuous $q$-Legendre polynomial is proved analytically by Rahman and Verma
\cite{20}, and a quantum $SU(2)$ group theoretic proof of this addition
formula is given by Koelink \cite{10}. However, the quantum
group theoretic proof more or less uses knowledge concerning
the structure of the addition formula for the continuous
$q$-Legendre polynomials. On the other hand, Koornwinder's \cite{13}
addition formula for the little $q$-Legendre polynomials follows
naturally from the interpretation of the little $q$-Jacobi
polynomials on the quantum $SU(2)$ group and this formula would
have been hard to guess without this interpretation. Rahman
\cite{19}, knowing what to prove,
has given an analytic proof of the addition formula
for the little $q$-Legendre polynomials.
As a follow-up to Koornwinder's \cite{14} paper, in which
he establishes an interpretation of a two-parameter family of
Askey-Wilson polynomials as zonal spherical elements on the
quantum $SU(2)$ group, abstract addition formulas, i.e. involving
non-commuting variables, have been given by Noumi and Mimachi
\cite{17} (see also \cite{18}) and by Koelink \cite{10}.
As a result of this approach there is a (degenerate) addition formula
for the two-parameter family of Askey-Wilson polynomials, cf.
\cite{17}, \cite{10}.

The group theoretic proof of the addition formula for the Legendre
polynomials starts with the spin $l$ ($l\in\Z$) representation
$t^l$ of the group $SU(2)$. The matrix elements $t^l_{n,m}$ are
known in terms of Jacobi polynomials and the matrix element $t^l_{0,0}$
is expressible in terms of the Legendre polynomial. Moreover, $t^l_{0,0}$
is the zonal spherical function with respect to the one-parameter
subgroup $K=S(U(1)\times U(1))$ of $SU(2)$, i.e.
$t^l_{0,0}(gk)=t^l_{0,0}(kg)=t^l_{0,0}(g)$ for all $g\in SU(2)$ and for
all $k\in K$. Using the homomorphism property we get
$$
t^l_{0,0}(gh)=\sum_k t^l_{0,k}(g)t^l_{k,0}(h), \qquad \forall g,h\in SU(2),
\tag{1.1}
$$
which yields the addition formula for the Legendre polynomials.
We can also view \thetag{1.1} as an expression for the unique
(up to a scalar) function $SU(2) \ni g\mapsto t^l_{0,0}(gh)$
in the span of the matrix elements $t^l_{n,m}$, which is left $K$-invariant
and right $hKh^{-1}$-invariant. It is this view of \thetag{1.1}
we adopt in this paper.

This view of \thetag{1.1} implies that we are not using the
comultiplication in the quantum group theoretic derivation of the
addition formula, in contrast with the quantum group theoretic
proofs of addition formulas mentioned. We start with a formula
relating the unique (up to a scalar) zonal spherical element,
which is left and right invariant with respect to different quantum
``subgroups'', to the matrix elements of the standard irreducible
unitary representations of the quantum $SU(2)$ group.
This formula is proved by Koornwinder in his paper
\cite{14} on zonal spherical elements on the quantum $SU(2)$ group.
In \cite{14} Koornwinder interpreted a two-parameter family of
Askey-Wilson polynomials as zonal spherical elements on the quantum
$SU(2)$ group. For a suitable choice of the parameters a quantum
group theoretic interpretation of the big $q$-Legendre polynomials
is obtained, which is a quantum group analogue of \thetag{1.1}.

This identity involves non-commuting variables, so we use a representation
to obtain an identity for operators acting on a Hilbert space. By
letting these operators act on suitable vectors of the Hilbert space
and taking inner products we obtain in a natural way an addition
formula for the big $q$-Legendre polynomial. The addition formula
involves Al-Salam--Carlitz polynomials, little $q$-Jacobi
polynomials and dual $q$-Krawtchouk polynomials.
The big $q$-Legendre polynomial corresponds to the term $t^l_{0,0}(gh)$
on the left hand side of \thetag{1.1} and the little
$q$-Jacobi polynomials, respectively the dual $q$-Krawtchouk
polynomials, correspond to $t^l_{0,k}(g)$, respectively
$t^l_{k,0}(h)$, in \thetag{1.1}. The Al-Salam--Carlitz
polynomials stem from the non-commutativity.

The dual $q$-Krawtchouk
polynomial tends to the Krawtchouk polynomial as $q\uparrow 1$ and the
Krawtchouk polynomial can be rewritten as a Jacobi polynomial, cf.
Koornwinder \cite{11, \S 2}, Nikiforov and Uvarov \cite{15,
\S\S 12, 22}. On the level of basic hypergeometric series we can rewrite
the dual $q$-Krawtchouk polynomial as a rational function resembling
a Jacobi polynomial of argument $z/(1+z)$, cf. \cite{9, p. 429}.
From the addition formula we obtain an expression for the product of a little
$q$-ultraspherical polynomal times a dual $q$-Krawtchouk polynomial as
a $q$-integral transformation of the big $q$-Legendre polynomials.
We show that a special case of this addition formula is related to a
special case of the addition formula for little $q$-Legendre
polynomials, cf. \cite{13}.

Although our initial relation is a
special case of the initial relation for Koornwinder's second addition
formula for $q$-ultraspherical polynomials, which he announced in
\cite{14, remark 5.4}, the addition formula for the big $q$-Legendre
polynomial proved here is not a special case of that second addition
formula. This is due
to the fact that we use an infinite dimensional $\ast$-representation
on our initial relation, whereas Koornwinder uses a one-dimensional
$\ast$-representation to obtain the $q$-Legendre case of his
addition formula for $q$-ultraspherical polynomials.

It should be noted that there is an abstract addition formula for
the big $q$-Legendre polynomial as a special case of the general
abstract addition formula mentioned before, cf. \cite{10},
\cite{17}. It is (at present) unknown whether it is
possible to derive an addition formula for the big $q$-Legendre polynomials
from the abstract addition formula. It might give an extension
of the result presented in this paper.

This paper is organised as follows. In sections 2 and 3 we recall the
necessary information on basic hypergeometric orthogonal polynomials
and on the quantum $SU(2)$ group. The main result is proved in \S 4.
Finally, in section 5 the limit $q\uparrow 1$ is considered. This
limit transition can be handled with the devices developed by
Van Assche and Koornwinder \cite{22} to prove that the addition
and product formula for the little $q$-Legendre polynomials tend to
the familiar addition and product formula for the Legendre polynomial.

\head 2. Preliminaries on basic hypergeometric
orthogonal polynomials\endhead

The notation for $q$-shifted factorials and basic hypergeometric series
is taken from the book \cite{7} by Gasper and Rahman.
We will assume $q\in (0,1)$.

The big $q$-Jacobi polynomials were introduced by Andrews and Askey
\cite{3, \S 3} and are defined by
$$
P_n(x;a,b,c,d;q) = {}_3 \varphi_2 \left( {{q^{-n}, abq^{n+1}, qax/c}\atop
{qa,-qad/c}};q,q\right).
\tag{2.1}
$$
The polynomial $P_n(x;1,1,c,d;q)$ is the big $q$-Legendre polynomial.

The monic big $q$-Jacobi polynomials $\hP_n$ with $a=0$, $b=0$, can be obtained
as a limit case of \thetag{2.1}.
First calculate the coefficient of $x^n$ in \thetag{2.1}
and next apply \cite{7, (3.2.3)} before taking $a\to 0$, $b\to 0$.
We find
$$
\aligned
\hP_n(x;0,0,c,d;q) =&\, d^n q^{{1\over 2}n(n-1)}\,
{}_2\varphi_1 \left( {{q^{-n}, c/x}\atop{0}};q,
-{{qx}\over{d}} \right), \\
\hP_n(x;0,0,c,d;q) =&\, (-c)^n q^{{1\over 2}n(n-1)}\,
{}_2\varphi_1 \left( {{q^{-n}, -d/x}\atop{0}};q,
{{qx}\over{c}} \right).
\endaligned
\tag{2.2}
$$
These polynomials satisfy the three-term recurrence relation
$$
\aligned
x\hP_n (x;0,0,c,d;q)= &\, \hP_{n+1}(x;0,0,c,d;q)
                   + q^n(c-d)\hP_n(x;0,0,c,d;q)\\
&\, +q^{n-1}cd(1-q^n) \hP_{n-1}(x;0,0,c,d;q).
\endaligned
\tag{2.3}
$$
Comparison of \thetag{2.3}
with the three-term recurrence relation for
the Al-Salam--Carlitz polynomials, cf. \cite{1, \S 4},
\cite{6, Ch.VI, \S 10}, shows that these monic big $q$-Jacobi polynomials
are Al-Salam--Carlitz polynomials with dilated argument,
$\hP_n(x;0,0,c,d;q)=c^n U_n^{(-d/c)}(x/c;q)$.
The orthogonality relations for the $\hP_n(\cdot;0,0,c,d;q)$ can be
phrased as
$$
\aligned
\int_{-d}^c &\bigl( \hP_n \hP_m\bigr) (x;0,0,c,d;q)\, (qx/c,-qx/d;q)_\infty
\, d_qx \\
=&\, \delta_{n,m} q^{{1\over 2}n(n-1)} (cd)^n (q;q)_n (1-q)c
(q,-d/c,-qc/d;q)_\infty .
\endaligned
\tag{2.4}
$$
Here the $q$-integral is defined by, cf. \cite{7, \S 1.11},
$$
\int_a^b f(x)\, d_qx = \int_0^b f(x)\, d_qx - \int_0^a f(x)\, d_qx,
\qquad \int_0^a f(x)\, d_qx = a(1-q)\sum_{k=0}^\infty f(aq^k)q^k .
$$

We will also need the little $q$-Jacobi polynomials $p_n(x;a,b;q)$,
cf. Andrews and Askey \cite{2, \S 3}, \cite{3, \S 3}.
The little $q$-Jacobi polynomials
are big $q$-Jacobi polynomials with $c=1$ and $d=0$ and normalised
such that the value at $0$ is $1$. Explicitly,
$$
p_n(x;a,b;q)={}_2\varphi_1 \left( {{q^{-n}, q^{n+1}ab}\atop{qa}};q,qx
\right) .
\tag{2.5}
$$

The last set of orthogonal polynomials needed is the set of dual
$q$-Krawtchouk polynomials, cf. \cite{21, \S 4},
which is a special case of the  $q$-Racah polynomials, cf.
\cite{5, \S 4}.
$$
R_n(q^{-x}-s^{-1}q^{x-N};s,N; q) = {}_3\varphi_2 \left(
{{q^{-n},q^{-x},-s^{-1}q^{x-N}}\atop{q^{-N},0}};q,q \right)
\tag{2.6}
$$
for $n\in\{ 0,\ldots,N\}$.

\head 3. Results on the quantum $SU(2)$ group\endhead

Let $q\in (0,1)$ be a fixed number. The unital $\ast$-algebra
$\A$ is generated by the elements $\a$ and $\g$ subject to
the relations
$$
\aligned
&\a\g = q\g\a , \quad \a\gs = q\gs \a , \quad \g\gs = \gs\g , \\
&\as\a +\g\gs = 1, \quad \a\as + q^2 \g\gs = 1.
\endaligned
\tag{3.1}
$$
For $q\uparrow 1$ the algebra can be identified with the algebra
of polynomials on the group $SU(2)$. The algebra $\A$ is actually
a Hopf $\ast$-algebra. See \cite{12}, \cite{16} for
references to the literature.

The irreducible unitary corepresentations of the Hopf $\ast$-algebra
$\A$ have been classified. For each dimension $2l+1$, $l\in {1\over 2}
\Zp$, there is precisely one such corepresentation, which we denote
by $t^l=\bigl( t^l_{n,m}\bigr)$, $n,m\in\{ -l,-l+1,\ldots,l\}$.
The matrix coefficients $t^l_{n,m}\in\A$
are explicitly known in terms of little $q$-Jacobi polynomials. For
our purposes it suffices to have
$$
\aligned
t^l_{0,m} &= d^l_m (\as )^m p_{l-m}(\g\gs ;q^{2m},q^{2m};q^2) (-q\gs)^m \\
t^l_{0,-m} &= d^l_m \g^m p_{l-m}(\g\gs ;q^{2m},q^{2m};q^2) \a^m
\endaligned
\tag{3.2}
$$
with
$$
d^l_m = {{q^{-m(l-m)}}\over{(q^2;q^2)_m}} \sqrt{ {{(q^2;q^2)_{l+m}}\over
{(q^2;q^2)_{l-m}}} }
$$
for $l\in\Zp$, $m=0,\ldots ,l$. See \cite{12}, \cite{16} for this
result as well as for references to the literature.

Next we recall a special case of Koornwinder's result \cite{14,
theorem~5.2} on general
spherical elements on the quantum $SU(2)$ group. The case we consider is
the case $\tau\to\infty$ of \cite{14, theorem~5.2}. Explicitly,
the following identity in $\A$ is valid;
$$
\sum_{m=-l}^l q^{-m/2} c^{l,\sigma}_m t^l_{0,m} =
C_l(\sigma) P_l(\r ;1,1,q^{2\sigma},1;q^2).
\tag{3.3}
$$
where $\sigma\in\R$,
$$
\align
c^{l,\sigma}_m =
c^{l,\sigma}_{-m} &= {{i^mq^{-(l+\sigma)m +{1\over 2}m^2}}\over{
\sqrt{(q^2;q^2)_{l+m}(q^2;q^2)_{l-m}}}} R_{l-m}(q^{-2l}-q^{-2l-2\sigma};
q^{2\sigma},2l;q^2), \\
C_l(\sigma) &= (-1)^l q^{-l^2-l} {{(-q^{2-2\sigma};q^2)_l}\over{
(q^{2l+2};q^2)_l}}
\endalign
$$
are constants and
$$
\r =\lim_{\tau\to\infty} 2q^{\sigma+\tau-1}\rst
= iq^{\sigma}(\as\gs -\g\a )-(1-q^{2\sigma})\gs\g \in\A .
$$
Here $\rst$ is defined in \cite{14, (4.8)}.
Equation \thetag{3.3} can be proved by redoing
Koornwinder's \cite{14} analysis with $X_\tau$ replaced
by $X_\infty$ or by taking the limit $\tau\to\infty$ in his
result \cite{14, theorem~5.2}. In the latter case we use the
limit transition of the Askey-Wilson polynomials to the big $q$-Legendre
polynomials as described in \cite{14, theorem~6.2},
$c^{l,\sigma}_l=i^lq^{-{1\over 2}l^2-l\sigma} (q^2;q^2)_{2l}^{-1/2}$, and
the limit
$$
\lim_{\tau\to\infty} q^{2\tau l}c^{l,\tau}_m =
{{(-1)^l \delta_{m,0}}\over{(q^{2l+2};q^2)_l}}.
$$
This follows for $m\geq 0$ from \cite{7, (3.2.3) with $e=0$, (1.5.3)}
and by the symmetry $c^{l,\tau}_{-m}=c^{l,\tau}_m$ for all $m$.

A $\ast$-representation $\pi$ of the commutation relations \thetag{3.1}
is acting on $\ell^2(\Zp )$ equipped with an orthonormal basis
$\{ e_n\}_{\{ n\in\Zp\} }$, and the explicit action of the generators
is given by
$$
\pi(\a )e_n = \sqrt{1-q^{2n}} e_{n-1}, \quad
\pi(\g )e_n = q^n e_n.
\tag{3.4}
$$
The irreducible $\ast$-representations of $\A$ have been classified, cf.
\cite{12} and the references therein. The infinite dimensional
$\ast$-representations are parametrised by the unit circle;
$\pi_\theta(\a)=\pi(\a)$ and
$\pi_\theta(\g )=e^{i\theta}\pi(\g )$ for $\theta\in [0,2\pi)$.

\head 4. Addition formula for big $q$-Legendre polynomials\endhead

In this section we prove an addition formula for the big $q$-Legendre
polynomials. We start by representing the relation \thetag{3.3} in $\A$
as an identity for operators in the Hilbert space $\ell^2(\Zp )$. Letting
these operators act on suitable vectors and taking inner products yields
the addition formula. This addition formula involves Al-Salam--Carlitz
polynomials, little $q$-Jacobi polynomials and dual $q$-Krawtchouk
polynomials. From the addition formula we find a $q$-integral representation
for the product of a little $q$-Jacobi polynomial and a dual $q$-Krawtchouk
polynomial.

Consider the action
of the infinite dimensional $\ast$-representation $\pi$ in $\ell^2(\Zp )$
on $\r$. The operator $\pi(\r )$ is a bounded self-adjoint operator and
the action on a basis vector $e_n$ of the standard orthonormal basis
is given by
$$
\pi(\r )e_n = -iq^{\sigma +n-1} \sqrt{1-q^{2n}} e_{n-1} -
q^{2n}(1-q^{2\sigma}) e_n + iq^{\sigma +n}\sqrt{1-q^{2n+2}}e_{n+1} .
$$
Consequently, $\sum_{n=0}^\infty p_n e_n$ is an eigenvector of $\pi(\r )$
for the eigenvalue $\lambda$ if and only if
$$
\lambda p_n = -iq^{\sigma +n} \sqrt{1-q^{2n+2}} p_{n+1} -
q^{2n}(1-q^{2\sigma}) p_n + iq^{\sigma +n-1}\sqrt{1-q^{2n}}p_{n-1}
\quad \forall\, n.
\tag{4.1}
$$
Since $p_{-1}=0$ and $p_0=1$, we view \thetag{4.1}
as a three-term recurrence for polynomials in $\lambda$. In order to determine
the polynomials from \thetag{4.1}
we calculate the leading coefficient $lc(p_n)=i^nq^{-\sigma n}
q^{-{1\over 2}n(n-1)} (q^2;q^2)_n^{-{1\over 2}}$ and determine the
three-term recurrence relation for the monic polynomials $\hp_n$;
$$
\lambda \hp_n(\lambda) = \hp_{n+1}(\lambda) -q^{2n}(1-q^{2\sigma})
\hp_n(\lambda)+(1-q^{2n})q^{2\sigma+2n-2}\hp_{n-1}(\lambda) .
\tag{4.2}
$$
Comparison of \thetag{4.2}
with the three-term recurrence relation \thetag{2.3}
for the big $q$-Jacobi polynomials with $a=0$ and $b=0$ leads to
$$
p_n(\lambda)= i^nq^{-\sigma n}q^{-{1\over 2}n(n-1)}
(q^2;q^2)_n^{-{1\over 2}} \hP_n(\lambda;0,0,q^{2\sigma},1;q^2) .
\tag{4.3}
$$
Denote the corresponding vector by $v_\lambda=\sum_{n=0}^\infty
p_n(\lambda)e_n$.

\proclaim{Proposition~4.1}
For $\lambda=-q^{2x}$, $x\in\Zp$, and
$\lambda = q^{2\sigma + 2x}$, $x\in\Zp$,
the vectors $v_\lambda$ constitute an orthogonal
basis of $\ell^2(\Zp )$.
\endproclaim

\demo{Proof} From the asymptotic formula, cf. \cite{8, (1.17)}, as
$n\to\infty$
$$
\hP_n(\lambda;0,0,c,d;q) \sim \lambda^n (c/\lambda,-d/\lambda ;q)_\infty
$$
for $\lambda\not=0$, $\lambda\not= cq^x$ and $\lambda\not= -dq^x$,
$x\in\Zp$, it follows that
$v_\lambda \not\in \ell^2(\Zp )$ for $\lambda\not=-q^{2x}$ and
$\lambda\not= q^{2\sigma+2x}$, $x\in\Zp$.

In the remaining cases we use the straightforward estimate
$$
\Bigl\vert {}_2\varphi_1 \left( {{q^{-n}, q^{-x}}\atop{0}};q,z
\right) \Bigr\vert \leq q^{-xn} (-q^{-x};q)_x (q,-\vert z\vert ;q)_\infty ,
\tag{4.4}
$$
for fixed $x\in\Zp$, in combination with the series representation
\thetag{2.2} for the
monic big $q$-Jacobi polynomials $\hP_n(\cdot;0,0,c,d;q)$
to see that we obtain eigenvectors $v_\lambda\in \ell^2(\Zp )$
for $\pi(\r )$ for the eigenvalues $\lambda=-q^{2x}$, $x\in\Zp$, and
$\lambda=q^{2\sigma+2x}$, $x\in\Zp$.

The orthogonality follows, since the vectors are eigenvectors of a
self-adjoint operator for different eigenvalues. It remains to prove
the completeness of the set of eigenvectors in $\ell^2(\Zp )$. To do
this we first calculate the length of the eigenvectors in $\ell^2(\Zp )$.
Consider $\lambda$ of the form $q^{2\sigma + 2x}$, $x\in\Zp$, then we have
proved the orthogonality relations
$$
\aligned
h_x \delta_{x,y}= \sum_{n=0}^\infty {{q^{n(n-1)} q^{-2\sigma n}}\over
{(q^2;q^2)_n}} \,& {}_2\varphi_1 \left( {{q^{-2x},q^{-2n}}\atop{0}};q^2,
-q^{2+2\sigma +2x}\right)\\
\times\,& {}_2\varphi_1 \left( {{q^{-2y},q^{-2n}}\atop{0}};q^2,
-q^{2+2\sigma +2y}\right) ,
\endaligned
\tag{4.5}
$$
for $x,y\in\Zp$, $h_x>0$.
We view the ${}_2\varphi_1$-series as a polynomial of degree $x$ in
the variable $q^{-2n}$. It has leading coefficient $(-1)^xq^{2x(x+\sigma)}$.
Since \thetag{4.5} holds, we have
$$
h_x = (-1)^xq^{2x(x+\sigma)}
\sum_{n=0}^\infty {{q^{n(n-1)} q^{-2\sigma n}}\over
{(q^2;q^2)_n}} \, {}_2\varphi_1 \left( {{q^{-2x},q^{-2n}}\atop{0}};q^2,
-q^{2+2\sigma +2x}\right) q^{-2nx}.
\tag{4.6}
$$
In \thetag{4.6} we replace the ${}_2\varphi_1$-series by
its terminating series representation
$$
\sum_{k=0}^x {{(q^{-2x};q^2)_k (q^{-2n};q^2)_k}\over{(q^2;q^2)_k}}
(-1)^k q^{2k(1+\sigma+x)}
$$
and we interchange the summations, which is justified by
the estimate \thetag{4.4}. The inner sum over $n$ starts at $n=k$ and
after a shift in the summation parameter the inner sum can be
evaluated using ${}_0\varphi_0(-;-;q,z)=(z;q)_\infty$,
cf. \cite{7, (1.3.16)}. The remaining
sum over $k$ can be summed using the $q$-binomial theorem
${}_1\varphi_0(q^{-p};-;q,z)=(q^{-p}z;q)_p$, cf. \cite{7, (1.3.14)}.
The result is
$$
h_x=q^{-2x} (q^2;q^2)_x (-q^{2\sigma+2};q^2)_x (-q^{-2\sigma};q^2)_\infty .
\tag{4.7}
$$
So
$w_x=v_{q^{2\sigma+2x}}/\parallel v_{q^{2\sigma+2x}}\parallel$
is an eigenvector of length $1$ of the self-adjoint operator $\pi(\r )$.

The orthogonality relations for the eigenvectors corresponding to
eigenvalues of the form $-q^{2x}$, $x\in\Zp$, is \thetag{4.5}
with $\sigma$ replaced by $-\sigma$. So
$u_x=v_{-q^{2x}}/\parallel v_{-q^{2x}}\parallel$
is an eigenvector of length $1$ of the self-adjoint operator $\pi(\r )$.

The set of orthonormal eigenvectors $\{ u_x\}_{\{ x\in\Zp\} } \cup
\{ w_x\}_{\{ x\in\Zp\} }$ forms a complete set of basis vectors for
$\ell^2(\Zp )$ if and only if the dual orthogonality relations
$$
\delta_{n,m} = \sum_{x=0}^\infty \langle u_x ,e_n\rangle
\overline{\langle u_x, e_m\rangle} \ +
\sum_{x=0}^\infty \langle w_x ,e_n\rangle
\overline{\langle w_x, e_m\rangle}
\tag{4.8}
$$
hold. It is easily seen that \thetag{4.8}
is equivalent to the orthogonality relations \thetag{2.4}
for the monic big $q$-Jacobi polynomials
$\hP_n(\cdot;0,0,q^{2\sigma},1;q^2)$. The first sum in \thetag{4.8}
corresponds to the $q$-integral over $[-1,0]$ and the second sum
corresponds to the $q$-integral over $[0,q^{2\sigma}]$. \qed
\enddemo

The orthogonality relations $\langle u_x,u_y\rangle =\delta_{x,y}$
(or $\langle w_x,w_y\rangle =\delta_{x,y}$) and
$\langle u_x,w_y\rangle=0$ can be stated in terms of
the $q$-Charlier polynomials, cf. \cite{7, exercise~7.13}.
(Note that the factor on the right hand side in
\cite{7, exercise~7.13} has to be replaced by its reciprocal.)

\proclaim{Corollary~4.2} Define the $q$-Charlier polynomials by
$$
c_n(x;a;q) = {}_2\varphi_1 ( q^{-n},x;0;q,-q^{n+1}/a), \qquad a> 0,
$$
then
$$
\sum_{x=0}^\infty {{a^x q^{{1\over 2}x(x-1)}}\over{(q;q)_x}}
\bigl( c_nc_m\bigr) (q^{-x};a;q) = \delta_{n,m}q^{-n} (q;q)_n (-q/a;q)_n
(-a;q)_\infty
$$
and
$$
\sum_{x=0}^\infty {{(-1)^x q^{{1\over 2}x(x-1)}}\over{(q;q)_x}}
c_n(q^{-x};a;q)c_m(q^{-x};a^{-1};q) = 0.
$$
\endproclaim

In order to convert \thetag{3.3} into a relation involving commuting
variables we apply the infinite dimensional $\ast$-representation $\pi$ to it.
We let the resulting bounded operator act on a standard basis vector $e_p$
and we take inner products with an eigenvector $v_\lambda\in\ell^2(\Zp)$,
cf. proposition~4.1. Next we use
the fact that $\pi$ is a $\ast$-representation to get the
following identity
$$
\sum_{m=-l}^l q^{-m/2} c^{l,\sigma}_m \langle
 \pi\bigl( t^l_{0,m}\bigr) e_p, v_\lambda\rangle
= C_l(\sigma)\P \langle e_p,v_\lambda\rangle,
\tag{4.9}
$$
since $\P$ is a polynomial with real coefficients and $\pi(\r)$ is
self-adjoint.
The operator on the left hand side of \thetag{4.9}
can be calculated explicitly by \thetag{3.2} and
\thetag{3.4}, since the standard basis vector $e_p$
is an eigenvector of $\pi(\g )$. Explicitly, for $m\geq 0$,
$$
\aligned
\pi\bigl(t^l_{0,m}\bigr) e_p& = d^l_m
(-1)^m q^{m(p+1)} \sqrt{(q^{2p+2};q^2)_m}
p_{l-m}(q^{2p};q^{2m},q^{2m};q^2) e_{p+m}, \\
\pi\bigl(t^l_{0,-m}\bigr) e_p& =d^l_m
q^{m(p-m)}  \sqrt{(q^{2p};q^{-2})_m}
p_{l-m}(q^{2(p-m)};q^{2m},q^{2m};q^2) e_{p-m},
\endaligned
\tag{4.10}
$$
with the convention $e_n=0$ for $n<0$.
Furthermore, from \thetag{4.3} it follows that for all $p\in\Zp$
$$
\langle e_p, v_\lambda \rangle = {{i^{-p} q^{-\sigma p}q^{-{1\over 2}p(p-1)}}
\over{\sqrt{(q^2;q^2)_p}}} \hP_p(\lambda;0,0,q^{2\sigma},1;q^2) .
\tag{4.11}
$$

Now we use \thetag{4.10} and \thetag{4.11} in
\thetag{4.9} together with the explicit values for
$C_l(\sigma)$, $c^{l,\sigma}_m$ and $d^l_m$, cf. \thetag{3.2},
\thetag{3.3}. Divide the resulting identity by
the factor in front of the monic big $q$-Jacobi polynomial in
\thetag{4.11} to obtain
$$
\aligned
&(-1)^lq^{-l^2-l}{{(-q^{2-2\sigma};q^2)_l}\over{(q^{2l+2};q^2)_l}} \P
\hP_l(\lambda;0,0,q^{2\sigma},1;q^2) =  \\
&{1\over{(q^2;q^2)_l}} R_l(q^{-2l}-q^{-2l-2\sigma};q^{2\sigma},2l;q^2)
p_l(q^{2p};1,1;q^2)\hP_p(\lambda;0,0,q^{2\sigma},1;q^2) \\
+&\sum_{m=1}^l (-1)^m{{q^{2m(p-l)}(q^{2p};q^{-2})_m}\over{(q^2;q^2)_{l-m}
(q^2;q^2)_m}} R_{l-m}(q^{-2l}-q^{-2l-2\sigma};q^{2\sigma},2l;q^2) \\
&\qquad\qquad\qquad\times p_{l-m}(q^{2(p-m)};q^{2m},q^{2m};q^2)
\hP_{p-m}(\lambda;0,0,q^{2\sigma},1;q^2) \\
+&\sum_{m=1}^l (-1)^m {{q^{m(m+1)-2m(\sigma+l)}}\over{(q^2;q^2)_{l-m}
(q^2;q^2)_m}} R_{l-m}(q^{-2l}-q^{-2l-2\sigma};q^{2\sigma},2l;q^2) \\
&\qquad\qquad\qquad\times p_{l-m}(q^{2p};q^{2m},q^{2m};q^2)
\hP_{p+m}(\lambda;0,0,q^{2\sigma},1;q^2)
\endaligned
\tag{4.12}
$$
for $\lambda=-q^{2x}$, $x\in\Zp$, or $\lambda=q^{2\sigma+2x}$, $x\in\Zp$.

We can now state and prove the main theorem of the paper.

\proclaim{Theorem~4.3}
{\rm (Addition formula for the big $q$-Legendre polynomial)}
With the notation of \thetag{2.1}, \thetag{2.2}, \thetag{2.5}
and \thetag{2.6} we have for $c,d>0$, $p,l\in\Zp$, $x\in\C$,
$$
\aligned
&(-1)^l q^{-{1\over 2}l(l+1)}{{(-qd/c;q)_l}\over{(q^{l+1};q)_l}}
P_l(x;1,1,c,d;q)\hP_p(x;0,0,c,d;q) = \\
&(q;q)_l^{-1} R_l(q^{-l}-{d\over c}q^{-l};{c\over d},2l;q)
p_l(q^p;1,1;q)\hP_p(x;0,0,c,d;q) \\
+& \sum_{m=1}^l (-1)^m{{d^mq^{m(p-l)}(q^p;q^{-1})_m}\over{(q;q)_{l-m}(q;q)_m}}
R_{l-m}(q^{-l}-{d\over c}q^{-l};{c\over d},2l;q) \\
&\qquad\qquad\qquad\times
p_{l-m}(q^{p-m};q^m,q^m;q)\hP_{p-m}(x;0,0,c,d;q) \\
+& \sum_{m=1}^l (-1)^m {{q^{{1\over
2}m(m+1)-lm}}\over{c^m(q;q)_{l-m}(q;q)_m}}
R_{l-m}(q^{-l}-{d\over c}q^{-l};{c\over d},2l;q) \\
&\qquad\qquad\qquad\times
p_{l-m}(q^p;q^m,q^m;q)\hP_{p+m}(x;0,0,c,d;q)
\endaligned
\tag{4.13}
$$
\endproclaim

\demo{Proof} Since \thetag{4.12}
only involves polynomials, it holds for all values of $\lambda$. In
\thetag{4.12} we replace $q^2$, $q^{2\sigma}$, $\lambda$ by
$q$, $c/d$, $x/d$. Now \thetag{4.13} follows from
$$
\aligned
P_n(x/d;a,b,c/d,1;q) &= P_n(x;a,b,c,d;q), \\
\hP_n(x/d;0,0,c/d,1;q) &= d^{-n} \hP_n (x;0,0,c,d;q),
\endaligned
$$
which is a consequence of \thetag{2.1} and \thetag{2.2}. \qed
\enddemo

\demo{Remark} 1. The choice of the infinite dimensional $\ast$-representation
does not influence the result.
We would obtain the same addition theorem if we had considered the
development of a $(\infty,\tau)$-spherical element in terms of
the standard matrix elements instead of \thetag{3.3}. \par
2. If we specialise $c=1$ and $d=0$ in \thetag{4.13}, then we can
sum the dual $q$-Krawtchouk polynomials $R_{l-m}$ by the $q$-Chu-Vandermonde
sum \cite{7, (1.5.3)}, from which we see that $R_{l-m}$ equals
$(q^{m+1};q)_{l-m}/(q^{l+m+1};q)_{l-m}$. The monic big
$q$-Jacobi polynomial $\hP_p$ with $a=b=d=0$, $c=1$ is summable
by the $q$-binomial theorem \cite{7, (1.3.14)}, which
results in $(-1)^pq^{{1\over 2}p(p-1)} (q^{1-p}x;q)_p$.
Furthermore, the big $q$-Legendre polynomial reduces to
$(-1)^lq^{{1\over 2}l(l+1)}p_l(x;1,1;q)$, so that we obtain the following
special case of \thetag{4.13};
$$
p_l(x;1,1;q)(q^{1-p}x;q)_p =
\sum_{m=0}^l {{q^{m(m-l+p)}(q;q)_{l+m}}\over{(q;q)_{l-m}(q;q)_m^2}}
p_{l-m}(q^p;q^m,q^m;q)(q^{1-p-m}x;q)_{p+m}.
$$
This corresponds to the case $x\to\infty$ of Koornwinder's addition
formula for the little $q$-Legendre polynomials \cite{13,
theorem~4.1 with $q^z=x$}.
\enddemo

The following $q$-integral representation for the product of a dual
$q$-Krawtchouk polynomial and a little $q$-ultraspherical
polynomial is a direct consequence of theorem~4.3 and the
orthogonality relations \thetag{2.4}. Just multipy
\thetag{4.13} by $\hP_{p+m}(x;0,0,c,d;q)$ and $q$-integrate
over $[-d,c]$ with respect to the weight function $(qx/c,-qx/d;q)_\infty$.

\proclaim{Corollary~4.4} For $c,d>0$, $p,l\in\Zp$,
$m\in\{ 0,\ldots,l\}$ we have
$$
\aligned
&R_{l-m}(q^{-l}-{d\over c}q^{-l};{c\over d},2l;q)
p_{l-m}(q^p;q^m,q^m;q) = \\
&C \int_{-d}^c P_l(x;1,1,c,d;q) \bigl( \hP_p\hP_{p+m}\bigr)(x;0,0,c,d;q)
(qx/c,-qx/d;q)_\infty d_qx
\endaligned
\tag{4.14}
$$
with
$$
C = {{(-1)^{l+m} q^{-{1\over 2}l(l+1) - {1\over 2}p(p-1)+m(l-p-m)}
c^{-p}d^{-p-m} (-qd/c;q)_l(q;q)_{l-m}}\over{(1-q)c(q^{l+1};q)_l (q^{m+1};q)_p
(q,-d/c, -qc/d;q)_\infty}}.
$$
\endproclaim

Multiplying \thetag{4.13} by $\hP_{p-m}(x;0,0,c,d,;q)$ and
$q$-integrating over $[-d,c]$ yields the same result \thetag{4.14}.
Specialising $m=0$ in \thetag{4.14} shows that the product of
the little $q$-Legendre polynomial and a dual $q$-Krawtchouk polynomial
can be written as a $q$-integral transform with a positive kernel of the
big $q$-Legendre polynomial.

\head 5. The limit case $q\uparrow 1$\endhead

In this section we show that the addition formula for the big
$q$-Legendre polynomials \thetag{4.13} and the product
formula \thetag{4.14} tend to the addition and
product formula for the Legendre polynomials as $q\uparrow 1$.
The general theorems of Van Assche and Koornwinder
\cite{22} used to obtain the addition and product formula
for the Legendre polynomials form the addition and product
formula for the little $q$-Legendre polynomials, cf. \cite{13},
are applicable in this case as well. See Askey \cite{4, Lecture~4} for
information on addition formulas for classical orthogonal polynomials.

We use the notation $R_n^{(\alpha,\beta)}(x)$ for the Jacobi
polynomial normalised by $R_n^{(\alpha,\beta)}(1)=1$.
First we note that the little and big $q$-Jacobi polynomials tend to the
Jacobi polynomials of shifted argument as $q\uparrow 1$;
$$
\aligned
\lim_{q\uparrow 1} P_n(x;q^\alpha,q^\beta,c,d;q) &=\,
R_n^{(\alpha,\beta)} \Bigl( {{2x+d-c}\over{c+d}} \Bigr),  \\
\lim_{q\uparrow 1} p_n(x;q^\alpha,q^\beta;q) &=\,
R_n^{(\alpha,\beta)} ( 1-2x).
\endaligned
\tag{5.1}
$$
The dual $q$-Krawtchouk polynomial can be rewritten as a
${}_2\varphi_2$-series, which tends to a Jacobi polynomial as
$q\uparrow 1$. This has also been used in \cite{9, p.~429}
to prove that the $q$-Krawtchouk polynomial tends to Jackson's
$q$-Bessel function. We can also let the dual $q$-Krawtchouk tend
to the Krawtchouk polynomial and use the relation between Krawtchouk
polynomials and Jacobi polynomials, cf. \cite{11, \S 2},
\cite{15, \S\S 12, 22}. The result is
$$
\lim_{q\uparrow 1} R_{l-m}(q^{-l}-{d\over c}q^{-l}; {c\over d}, 2l;q)
= {{(m+1)_{l-m}}\over{(l+m+1)_{l-m}}}(1+{d\over c})^{l-m}
R_{l-m}^{(m,m)}\Bigl( {{c-d}\over{c+d}} \Bigr).
\tag{5.2}
$$

In order to apply the theorems of Van Assche and Koornwinder \cite{22}
we have to consider the orthonormal big $q$-Jacobi polynomials with $a=0$,
$b=0$. Define
$$
p_k(x;q)= {{\hP_k(x;0,0,c,d;q)}\over{ q^{{1\over 4}k(k-1)}
(cd)^{k/2}(1-q)^{1/2} c^{1/2}\sqrt{(q;q)_k (q,-d/c,-qc/d;q)_\infty}}},
\tag{5.3}
$$
then the polynomials $p_k(x;q)$ satisfy the recurrence relation
$$
xp_k(x;q)=a_{k+1}(q) p_{k+1}(x;q) + b_k(q) p_k(x;q) + a_k(q) p_{k-1}(x;q)
$$
with
$$
a_k(q) = q^{{1\over 2}(k-1)}\sqrt{cd(1-q^k)}, \qquad
b_k(q) = q^k(c-d).
$$
Fix $r\in (0,1)$ and define $a_{k,n}=a_k(r^{1/n})$ and $b_{k,n}=b_k(r^{1/n})$.
The following limits are easily established;
$$
\lim_{n\to\infty} a_{n,n}=\sqrt{rcd(1-r)}>0,\qquad\quad
\lim_{n\to\infty} b_{n,n}=r(c-d)\in\R
$$
and
$$
\lim_{n\to\infty} (a_{k,n}^2-a^2_{k-1,n})=0,\qquad\quad
\lim_{n\to\infty} (b_{k,n}-b_{k-1,n})=0
$$
uniformly in $k$. Now \cite{22, theorem~1} can be applied and it yields
$$
\lim_{n\to\infty} {{p_{n+1}(x;r^{1/n})}\over{p_n(x;r^{1/n})}} =
\rho\Bigl( {{x-r(c-d)}\over{2\sqrt{rcd(1-r)}}}\Bigr)
\tag{5.4}
$$
uniformly on compact subsets of $\C\backslash [-d,c]$. Here
$\rho(x)=x+\sqrt{x^2-1}$ and the square root is the one for which
$\vert \rho(x)\vert>1$ for $x\notin [-1,1]$.
Rewriting \thetag{5.4} in terms of the big $q$-Jacobi polynomial and
iterating yields
$$
\lim_{p\to\infty} {{\hP_{p+m}(x;0,0,c,d;r^{1/p})}\over{
\hP_p(x;0,0,c,d;r^{1/p})}} = \bigl( cdr(1-r)\bigr)^{m/2}
\rho^m\Bigl( {{x-r(c-d)}\over{2\sqrt{rcd(1-r)}}}\Bigr)
\tag{5.5}
$$
for all $m\in\Z$ and $x\in \C\backslash [-d,c]$.

Now the proof that \thetag{4.13} tends to the addition formula for
Legendre polynomials can be finished. Replace $q$ by $r^{1/p}$ in
\thetag{4.13}, divide both sides by
$(q^{l+1};q)_l^{-1}\hP_p(x;0,0,c,d;r^{1/p})$ and let $p\to\infty$, i.e.
$q\uparrow 1$, then we can use \thetag{5.1}, \thetag{5.2}
and \thetag{5.5} to obtain, after a short calculation,
$$
\aligned
&\qquad R_l^{(0,0)}\Bigl( {{2x+d-c}\over{c+d}}\Bigr) =
(-1)^l R_l^{(0,0)}\Bigl( {{c-d}\over{c+d}}\Bigr) R^{(0,0)}_l(1-2r) \\
+&\sum_{m=1}^l{{(l+1)_l(m+1)_{l-m}}\over{(l-m)!m!(l+m+1)_{l-m}}}
(-1)^{l+m}(1+{d\over c})^{-m} \bigl({d\over c}r(1-r)\bigr)^{m/2} \\
& \times R^{(m,m)}_{l-m}\Bigl( {{c-d}\over{c+d}}\Bigr)
R^{(m,m)}_{l-m}(1-2r) \Bigl[
\rho^m\Bigl( {{x-r(c-d)}\over{2\sqrt{rcd(1-r)}}}\Bigr) +
\rho^{-m}\Bigl( {{x-r(c-d)}\over{2\sqrt{rcd(1-r)}}}\Bigr)\Bigr].
\endaligned
\tag{5.6}
$$
The term in square brackets equals
$2T_m\bigl((x-r(c-d))/2\sqrt{rcd(1-r)}\bigr)$,
where $T_m(\cos\theta)=\cos m\theta$ is the Chebyschev polynomial
of the first kind. In \thetag{5.6} we also use
$R^{(m,m)}_n(-x)=(-1)^nR^{(m,m)}_n(x)$, then we find, after a short
manipulation of the Pochhammer symbols,
$$
\aligned
&R_l^{(0,0)}\Bigl( {{2x+d-c}\over{c+d}}\Bigr) =
R_l^{(0,0)}\Bigl( {{c-d}\over{c+d}}\Bigr) R^{(0,0)}_l(1-2r) \\
+&2\sum_{m=1}^l {{(l+m)!}\over{(l-m)!(m!)^2}}
(1+{d\over c})^{-m} \bigl({d\over c}r(1-r)\bigr)^{m/2}
R^{(m,m)}_{l-m}\Bigl( {{d-c}\over{c+d}}\Bigr) \\
&\qquad\qquad\times R^{(m,m)}_{l-m}(1-2r) T_m\Bigl(
{{x-r(c-d)}\over{2\sqrt{rcd(1-r)}}}\Bigr).
\endaligned
\tag{5.7}
$$
Since the dependence on $x$ in \thetag{5.7} is polynomial, the
restriction $x\in\C\backslash [-d,c]$ can be removed.
Formula \thetag{5.7} is equivalent to the addition formula for the
Legendre polynomial, cf. \cite{4, Lecture~4},
$$
\aligned
&\qquad R_l^{(0,0)}\bigl( xy+t\sqrt{(1-x^2)(1-y^2)}\bigr)
= R^{(0,0)}_l(x)R^{(0,0)}_l(y) \\
+&2\sum_{m=1}^l {{(l+m)!}\over{(l-m)!(m!)^2}}2^{-2m}
\biggl(\sqrt{(1-x^2)(1-y^2)}\biggr)^m R^{(m,m)}_{l-m}(x)
R^{(m,m)}_{l-m}(y)T_m(t)
\endaligned
\tag{5.8}
$$
by identifying $(d-c)/(c+d)$, $1-2r$, $\bigl( x-r(c-d)\bigr)/2\sqrt{rcd(1-r)}$
with $x$, $y$ and $t$.

The limit case of the product formula \thetag{4.14} can also be
handled with the methods developed by Van Assche and Koornwinder
\cite{22}. Note that
$$
A=\lim_{n\to\infty} a_{n+k,n}= \sqrt{rcd(1-r)},\qquad
B=\lim_{n\to\infty} b_{n+k,n}= r(c-d)
$$
for all $k\in\Z$. Now \cite{22, theorem~2} can be applied
to yield
$$
\lim_{p\to\infty} \int_{-d}^c f(z)p_p(z;r^{1/p})p_{p+m}(z;r^{1/p})\,
d\mu_p(z) = {1\over\pi} \int_{B-2A}^{B+2A}
{{f(z) T_m\bigl( (z-B)/2A\bigr)}\over{\sqrt{4A^2-(z-B)^2}}}\, dz
\tag{5.9}
$$
for all continuous functions $f$ on $[-d,c]$. Here
$$
\int_{-d}^c f(z)d\mu_p(z)= \int_{-d}^c f(z) (qz/c,-qz/d;q)_\infty\, d_qz
$$
with $q$ on the right hand side replaced by $r^{1/p}$, and the $p_p(z;q)$
are the orthonormal big $q$-Jacobi polynomials with $a=b=0$,
cf. \thetag{5.3}.

If we now use \thetag{5.9}, \thetag{5.3}, \thetag{5.2}
and \thetag{5.1} to take the limit $q=r^{1/p}\uparrow 1$, i.e.
$p\to\infty$, in \thetag{4.14}, we obtain
$$
\align
&(-1)^{l+m} \bigl( dr(1-r)/c\bigr)^{-m/2} {{(l-m)!m!}\over{(l+1)_l\pi}}
\int_{B-2A}^{B+2A} R^{(0,0)}_l\Bigl({{2z+d-c}\over{c+d}} \Bigr)
{{T_m\bigl( (z-B)/2A\bigr)}\over{\sqrt{4A^2-(z-B)^2}}}\, dz \\
&\qquad ={{(m+1)_{l-m}}\over{(l+m+1)_{l-m}}}(1+{d\over c})^{-m}
R_{l-m}^{(m,m)}\Bigl( {{c-d}\over{c+d}} \Bigr)
R_{l-m}^{(m,m)}(1-2r)
\endalign
$$
with $A=\sqrt{rcd(1-r)}$, $B=r(c-d)$. By changing the integration variable
to $t=(z-B)/2A$, replacing $(d-c)/(c+d)$, $1-2r$ by $x$, $y$ and using
$R^{(m,m)}_n(-x)=(-1)^nR^{(m,m)}_n(x)$ we obtain the product formulas
$$
\aligned
R^{(m,m)}_{l-m}(x)R^{(m,m)}_{l-m}(y)=&  2^{2m}
{{(l-m)!(m!)^2}\over{\pi (l+m)!}} \bigl( \sqrt{(1-x^2)(1-y^2)}\bigr)^{-m} \\
&\times \int_{-1}^1 R^{(0,0)}_l\bigl( xy+t\sqrt{(1-x^2)(1-y^2)}\bigr)
{{T_m(t)}\over{\sqrt{1-t^2}}}\, dt.
\endaligned
\tag{5.10}
$$



\Refs

\ref\no 1
\by W.A.~Al-Salam and L.~Carlitz
\paper Some orthogonal $q$-polynomials
\jour Math. Nachr.
\vol 30
\yr 1965
\pages 47--61
\endref

\ref\no 2
\by G.E.~Andrews and R.~Askey
\paper Enumeration of partitions: The role of Eulerian series and
$q$-orthogonal polynomials
\inbook  in ``Higher Combinatorics''
\ed M.~Aigner
\publaddr Reidel
\yr 1977
\pages 3--26
\endref

\ref\no 3
\bysame 
\paper Classical orthogonal polynomials
\inbook  in ``Polyn\^omes Orthogonaux et Applications''
\bookinfo Lecture Notes Math. 1171
\eds C.~Brezinski, A.~Draux, A.P.~Magnus, P.~Maroni and A.~Ronveaux
\publaddr Springer
\yr 1985
\pages 36--62
\endref

\ref\no 4
\by R.~Askey
\book Orthogonal Polynomials and Special Functions
\bookinfo CBMS-NSF Regional Conference Series in Applied Mathematics {\bf 21}
\publ SIAM
\yr 1975
\endref

\ref\no 5
\by R.~Askey and J.~Wilson
\paper A set of orthogonal polynomials that generalize the Racah
coefficients or $6-j$ symbols
\jour SIAM J. Math. Anal.
\vol 10
\yr 1979
\pages 1008--1016
\endref

\ref\no 6
\by T.S.~Chihara
\book An Introduction to Orthogonal Polynomials
\bookinfo Mathematics and its Applications 13
\publaddr Gordon and Breach
\yr 1978
\endref

\ref\no 7
\by G.~Gasper and M.~Rahman
\book Basic Hypergeometric Series
\bookinfo Encyclopedia of Mathematics and its Applications 35
\publaddr Cambridge University Press
\yr 1990
\endref

\ref\no 8
\by M.E.H.~Ismail and J.A.~Wilson
\paper Asymptotic and generating relations for the $q$-Jacobi and
${}_4 \varphi_3$ polynomials
\jour J. Approx. Theory
\vol 36
\yr 1982
\pages 43--54
\endref

\ref\no 9
\by H.T.~Koelink
\paper Hansen-Lommel orthogonality
relations for Jackson's $q$-Bessel functions
\jour J. Math. Anal. Appl.
\vol 175
\yr 1993
\pages 425--437
\endref

\ref\no 10
\bysame 
\paper The addition formula for
continuous $q$-Legendre polynomials and associated spherical
elements on the $SU(2)$ quantum group related to Askey-Wilson
polynomials
\jour SIAM J. Math. Anal.
\vol 25
\yr 1994
\pages 197--217
\endref

\ref\no 11
\by T.H.~Koornwinder
\paper Krawtchouk polynomials, a unification of two different group theoretic
interpretations
\jour SIAM J. Math. Anal.
\vol 13
\yr 1982
\pages 1011--1023
\endref

\ref\no 12
\bysame 
\paper Orthogonal polynomials in connection with quantum groups
\inbook  in ``Orthogonal Polynomials: Theory and Practice''
\bookinfo NATO ASI series C, vol 294
\ed P.~Nevai
\publaddr Kluwer
\yr 1990
\pages 257--292
\endref

\ref\no 13
\bysame 
\paper
The addition formula for
little $q$-Legendre polynomials and the $SU(2)$ quantum group
\jour SIAM J. Math. Anal.
\vol 22
\yr 1991
\pages 195--301
\endref

\ref\no 14
\bysame 
\paper Askey-Wilson polynomials
as zonal spherical functions on the $SU(2)$ quantum group
\jour SIAM J. Math. Anal.
\vol 24
\yr 1993
\pages 795--813
\endref

\ref\no 15
\by A.F.~Nikiforov and V.B.~Uvarov
\book Special Functions of Mathematical Physics
\bookinfo Translated from the Russian by R.P.~Boas
\publaddr Birkh\"auser
\yr 1988
\endref

\ref\no 16
\by M.~Noumi
\paper Quantum groups and $q$-orthogonal polynomials.
Towards a realization of Askey-Wilson polynomials on $SU_q(2)$
\inbook  in ``Special Functions''
\bookinfo ICM-90 Satellite Conference Proceedings
\eds M.~Kashiwara and T.~Miwa
\publaddr Springer
\yr 1991
\pages 260--288
\endref

\ref\no 17
\by M.~Noumi and K.~Mimachi
\paper Askey-Wilson polynomials and the quantum group $SU_q(2)$
\jour Proc. Japan Acad., Ser. A
\vol 66
\yr 1990
\pages 146--149
\endref

\ref\no 18
\bysame 
\paper Askey-Wilson polynomials as spherical functions on $SU_q(2)$
\inbook  in ``Quantum Groups''
\bookinfo Lecture Notes Math. 1510
\ed P.P.~Kulish
\publaddr Springer
\yr 1992
\pages 98--103
\endref

\ref\no 19
\by M.~Rahman
\paper A simple proof of Koornwinder's addition formula for the
little $q$-Legendre polynomials
\jour Proc. Amer. Math. Soc.
\vol 107
\yr 1989
\pages 373--381
\endref

\ref\no 20
\by M.~Rahman and A.~Verma
\paper Product and addition
formulas for the continuous q-ultraspherical polynomials
\jour SIAM J. Math. Anal.
\vol 17
\yr 1986
\pages 1461--1474
\endref

\ref\no 21
\by D.~Stanton
\paper Orthogonal polynomials and Chevalley groups
\inbook  in ``Special Functions: Group Theoretical Aspects and Applications''
\eds R.A.~Askey, T.H.~Koornwinder and W.~Schempp
\publaddr Reidel
\yr 1984
\pages 87--128
\endref

\ref\no 22
\by W.~Van Assche and T.H.~Koornwinder
\paper Asymptotic behaviour for Wall polynomials and the addition formula
for little $q$-Legendre polynomials
\jour SIAM J. Math. Anal.
\vol 22
\yr 1991
\pages 302--311
\endref

\endRefs
\enddocument